\documentclass[11pt]{article}
\usepackage{amscd,amssymb,stmaryrd}
\usepackage{amsmath}
\usepackage{bbm}
\usepackage{graphicx}
\usepackage{subcaption}
\usepackage[utf8]{inputenc}
\usepackage[export]{adjustbox}
\usepackage{wrapfig}
\usepackage{amstext}
\usepackage{pstricks,pst-node,pst-plot,pst-coil}
\usepackage{amsthm}
\usepackage{amsmath}
\usepackage{color}
\usepackage{mathtools}
\usepackage{hyperref}
\usepackage[english]{babel}
\usepackage{booktabs}
\usepackage{mathrsfs}
\usepackage{comment}
\usepackage{float}
\usepackage[linesnumbered,ruled,vlined, ruled,vlined]{algorithm2e}
\usepackage[font=small,labelfont=bf]{caption}

\theoremstyle{definition}

\theoremstyle{remark}
\newtheorem*{remark}{Remark}

\def\eg{\emph{e.g., }}
\def\ie{\emph{i.e., }}
\def\etal{\emph{et al.}}

\title{DeepONet-Grid-UQ: A Trustworthy Deep Operator Framework for Predicting the Power Grid's Post-Fault Trajectories}
\author{Christian Moya, Shiqi Zhang, Meng Yue, and Guang Lin}

\begin{document}
\maketitle
\begin{abstract}
This paper proposes a new data-driven method for the reliable prediction of power system post-fault trajectories. The proposed method is based on the fundamentally new concept of Deep Operator Networks~(DeepONets). Compared to traditional neural networks that learn to approximate functions, DeepONets are designed to approximate nonlinear operators. Under this operator framework, we design a DeepONet to (1) take as inputs the fault-on trajectories collected, for example, via simulation or phasor measurement units, and (2) provide as outputs the predicted post-fault trajectories. In addition, we endow our method with a much-needed ability to balance efficiency with reliable/trustworthy predictions via uncertainty quantification. To this end, we propose and compare two methods that enable quantifying the predictive uncertainty. First, we propose a \textit{Bayesian DeepONet}~(B-DeepONet) that uses stochastic gradient Hamiltonian Monte-Carlo to sample from the posterior distribution of the DeepONet parameters. Then, we propose a \textit{Probabilistic DeepONet}~({\color{black}Prob-DeepONet}) that uses a probabilistic training strategy to equip DeepONets with a form of automated uncertainty quantification, at virtually no extra computational cost. Finally, we validate the predictive power and uncertainty quantification capability of the proposed B-DeepONet and {\color{black}Prob-DeepONet} using the IEEE 16-machine 68-bus system.
\end{abstract}

\section{Introduction} \label{sec:introduction}
Many critical infrastructure systems depend heavily on a robust, reliable, and always-on power grid. However, our power grid is constantly exposed to rare but severe events and disturbances that could compromise its operation and security. In the worst-case scenario, such disturbances could lead to instabilities, which, in turn, could trigger system-wide blackouts.  

To assess the power grid's dynamic security~\cite{kundur2007power}, system operators implement an offline transient analysis procedure known in practice as the dynamic security analysis (DSA), usually based on the $N$-1 criteria~\cite{alvarado2002transmission}. The $N$-1 based DSA (1) seeks to predict whether the power grid will remain safely operating after facing a significant disturbance (\eg a transmission line fault or a generator trip) selected from a set of credible disturbances and (2) requires simulating the power grid's dynamic response.  

To simulate the power grid's dynamic response, one must simulate a set of nonlinear differential-algebraic equations~(DAE)~\cite{kundur2007power}. Simulating this set of DAEs is, however, a very challenging task. Indeed, the classical explicit integration methods fail catastrophically on such a task~\cite{iserles2009first}. As a result, most of the software tools used for transient stability implement numerically stable integration methods~\cite{kundur2007power} designed specifically for DAEs. However, the computational cost and memory required for these methods are very high and constitute the main obstacle towards deploying a more efficient and possibly online transient analysis~\cite{schainker2006real}. With the transformation the power grid now faces (\eg by deploying a large number of distributed renewable energy resources or a liberalized market), soon, it will become imperative for electric utilities to assess its dynamic security near real-time, which calls for more efficient transient analysis methods.

\subsection{Prior works} \label{sub-sec:prior-works} 
Of course, to accelerate the simulation of DAEs, we can resort to schemes aiming at parallelizing the existing numerical methods. Parallel computing, however, consumes a high amount of resources, and the performance can quickly reach the saturation, which could prevent its application in near-real-time transient analysis. Besides parallel computing, other methods have been proposed to assess transient stability efficiently in the last decades. Here, we classify them into (1) transient analysis methods that exploit the mathematical structure of the power grid models and (2) Machine Learning~(ML)-enabled transient analysis methods.

\subsubsection{Power grid model-based transient analysis methods} An alternate approach to the traditional offline transient analysis method was proposed in the early '80s under the name of \textit{direct methods}~\cite{pai1981transient,varaiya1985direct,hiskens1989energy} (or energy function-based methods). These methods provide transient stability certification without the time-consuming numerical simulation of the post-fault dynamics. More specifically, direct methods assess transient stability by analyzing the power grid dynamical model with the aid of an energy function (\ie a Lyapunov-like function), which certifies convergence of the states to the stable equilibrium point. Modern versions of direct methods~\cite{chiang1994bcu,chiang2011direct} have been successfully engineered to the point that they were implemented at the utility level. However, the scalability and conservativeness of the classical energy-function method limit their applicability to even relatively large power grids.

Other works have also addressed the transient analysis problem using power grid dynamical models. For example, Dorfler \etal~\cite{dorfler2012synchronization} studied the stability of power grids with overdamped generators via studying the synchronization of Kuramoto oscillators. The authors of~\cite{jin2005power} used Hamilton Jacobi reachability to find an inner approximation of the region of attraction, which enables the online transient stability certification of initial conditions. Finally, \cite{anghel2013stability} proposed using the sum of squares programming, and \cite{caliskan2014compositional} proposed using a port Hamiltonian approach. All the above methods use simplified power grid dynamics. So, their applicability to more complex power grids remains unclear. In addition, most of these works produce only stability margins but not the power grid post-fault trajectories, which system operators and planners often need.  

\subsubsection{Machine Learning~(ML)-enabled transient analysis methods}
Machine Learning holds the promise to revolutionize science and engineering, including power engineering. Thus, as expected, many works have proposed to address the transient analysis problem using Machine Learning and Artificial Intelligence-inspired methodologies. The main idea behind most of these works is to learn a mapping from initial states/responses, collected via simulation or phasor measurement units~(PMU), to a certificate of transient stability. 

For example, He \etal~\cite{he2013robust} proposed using ensemble tree learning to construct a hierarchical characterization of the dangerous regions in the space of possible contingencies and operating/initial states. In~\cite{wehenkel1989artificial}, the authors proposed a classification approach to learn the decision rules for certifying transient stability from initial responses. Using recurrent neural networks, the authors in \cite{james2017intelligent} proposed a self-adaptive transient analysis method. In~\cite{gupta2018online}, the authors developed a convolutional neural network that takes as the input measurements collected using PMUs and outputs a transient stability certificate. However, the main drawback of these aforementioned ML-based approaches is that they focus primarily on developing a binary/discrete indicator that certifies whether the power grid is transient stable or not after a disturbance. 

To provide a more quantitative metric for transient stability, Zhu \etal~\cite{zhu2019hierarchical} developed a convolutional neural network that takes as the input the fault-on trajectory collected using PMUs and outputs not only a transient stability certificate but also an estimate of the stability margin. This ML-based tool can be fast, but it may be insufficient for system operators and planners. These operators and planners may be interested in knowing the trajectories of various states of the power grid after the disturbance. These trajectories are crucial to estimate whether any voltage or frequency will violate pre-defined limits and trigger, for example, load shedding. Hence, in the authors' previous work~\cite{li2020machine}, we developed a deep-learning framework that uses the long-short-term memory ~(LSTM) to predict the transient responses of power grid states. The effectiveness and speed of this approach in~\cite{li2020machine} was undoubtedly a first step towards the development of a more complete and efficient ML-based transient analysis tool for power grids. However, our previous work, like most of the proposed ML-based transient analysis works, (1) requires a considerable amount of supervision for successful deployment and (2) fails to tackle the inherent trade-off between being efficient and being reliable. 

\subsection{Our work}   \label{sub-sec:our-work}
This paper proposes a machine learning-based method for predicting the power grid's post-fault transient response to faults. However, compared to existing ML-based methods, we employ a different class of neural networks, introduced in the seminal work of Lu~\etal~\cite{lu2021learning}, that we train to approximate the operator arising from such transient response. Under this operator framework, our objectives are:
\begin{enumerate}
    \item \textit{Post-fault trajectory prediction:} deriving a Deep Operator Network~(DeepONet) framework that observes the fault-on trajectory (of a given state) to predict the post-fault trajectory (of the same state).
    \item \textit{Trustworthy prediction:} designing an effective strategy that attaches a measure of confidence/reliability to the predicted trajectory.
\end{enumerate}
Compared to~\cite{li2020machine}, we do not focus on improved mean trajectory prediction. One can always improve the mean predictive performance by using novel neural network architectures or enriching the datasets (which may not be possible in practice). Hence, in this paper, we consider a smaller dataset than in~\cite{li2020machine}. Furthermore, during training, we will not use information from neighboring buses, which generally improves the mean prediction. The above two assumptions make the proposed problem of providing trustworthy/reliable post-fault trajectory predictions much more challenging. Let us conclude this section by enumerating our contributions and describing the paper's organization.  
\begin{enumerate}
    \item We first design~(in Section~\ref{sec:operator-regression}) a data-driven deep operator network~(DeepONet) framework for post-fault trajectory prediction. Our DeepONet learns to approximate the solution operator that maps trajectories containing fault-on information to their corresponding post-fault trajectories. Once trained, the proposed DeepONet provides a more efficient way of calculating the post-fault trajectories than traditional transient response simulators. 
    \item We then endow the proposed DeepONet with the ability to provide trustworthy predictions. To this end, we develop two methods that quantify the heteroscedastic aleatoric uncertainty (\ie the input and data-inherent uncertainty) of the predicted trajectories. 
    \item In the first method, we design~(in Section~\ref{sec:bayesian-deeponet}) a Bayesian framework for DeepONets. Here the Bayesian DeepONet~(B-DeepONet) represents the prior of the trainable parameters. We then employ the (theoretically sound) Stochastic Gradient Hamiltonian Monte-Carlo~(SGHMC) algorithm to sample the corresponding posterior. We use the collected samples to quantify the uncertainty of predicted post-fault trajectories.
    \item For our second method, we propose~(in Section~\ref{sec:probabilistic-deeponet}) a probabilistic DeepONet~(Prob-DeepONet) framework. This framework uses a special DeepONet architecture and probabilistic training to transform the DeepONet predictions into point-wise Gaussian estimates. Predicting Gaussian estimates enables quantifying the uncertainty of the predicted trajectories at virtually no extra computational cost. 
    \item Finally, we test~(in Section~\ref{sec:numerical-experiments}) the trustworthiness and predictive power of the proposed B-DeepONet and Prob-DeepONet using the IEEE-{\color{black}68}-bus power grid. 
\end{enumerate}
We organize the rest of this paper as follows. In Section~\ref{sec:operator-regression}, we formulate the problem of predicting the post-fault trajectories as an operator regression problem. Moreover, we introduce the DeepONet to approximate such an operator. The B-DeepONet and Prob-DeepONet frameworks that enable trustworthy post-fault trajectory prediction are detailed, respectively, in Section~\ref{sec:bayesian-deeponet} and Section~\ref{sec:probabilistic-deeponet}. In Section~\ref{sec:numerical-experiments}, we test, using the IEEE-69 bus power grid, the reliability and predictive power of the proposed frameworks. Section~\ref{sec:discussion} discusses our results and future work.  Finally, Section~\ref{sec:conclusion} concludes the paper. 

\section{Operator Regression for Transient Response Prediction} \label{sec:operator-regression}
This section describes how to pose the post-fault trajectory prediction problem as an operator regression problem.

\textit{Post-fault trajectory prediction:} Faults and disturbances on transmission lines (or other large grid equipment) can drive the power grid to instability. In a typical fault scenario, the power grid experiences the fault and moves away from the pre-fault stable equilibrium until the system clears the fault~\cite{chiang2011direct}. Once the fault is cleared, the power grid experiences a transient response known as the post-fault dynamics. For power grid operators, it is then essential to know whether these post-fault dynamics converge or not to a stable equilibrium point~\cite{kundur2007power}. If the trajectories fail to converge, a series of conventional protection schemes will be triggered, leading, for example, to load shedding.

\textit{Operator regression:} We formalize next how to predict post-fault trajectories using operator regression. Assume the power grid was operating at a stable equilibrium point. Then, at some time~$t_f$, the grid experiences a fault, which results in a structural change due to the system protective actions~\cite{chiang2011direct}. Assume the resulting fault-on dynamics are confined to the time interval~$[t_f,t_{cl}]$, where $t_{cl}$ denotes the time when the fault is cleared. During~$[t_f, t_{cl}]$, we can describe the fault-on dynamics using the set of differential-algebraic equations~(DAE):
\begin{subequations} \label{eq:DAE-fault-on}
\begin{align}
\dot{x}(t) &= f_F(x(t),z(t)) \quad t \in [t_f, t_{cl}], \\
0 &= g_F(x(t),z(t)),    
\end{align}
\end{subequations}
where $x$ is the vector of dynamic states and $z$ the vector of algebraic states. Now, assume the system clears the fault at time~$t_{cl}$ and no additional protective actions occur after~$t_{cl}$. Then, the post-fault dynamics, \ie the dynamics within the interval $(t_{cl}, T]$\footnote{$T$ denotes a time (in the order of seconds) that can capture the power grid's transient response.}, are given by the DAEs:
\begin{subequations} \label{eq:DAE-post-fault}
\begin{align}
\dot{x}(t) &= f_{PF}(x(t),z(t)) \quad t \in (t_{cl},T], \\
0 &= g_{PF}(x(t),z(t)).    
\end{align}
\end{subequations}
In this paper, we aim to learn the operator, denoted as~$G^\dagger$, mapping a fault-on trajectory (of a selected dynamic or algebraic state~$w \in \{x,z\}$), which satisfies~\eqref{eq:DAE-fault-on}, to the post-fault trajectory (of the \textit{same} state~$w$), which satisfies~\eqref{eq:DAE-post-fault}, \ie 
\begin{align} \label{eq:operator}
    G^\dagger: w|_{[t_f,t_{cl}]} \mapsto w|_{(t_{cl},T]}. 
\end{align}
Here $w|_{[t_f,t_{cl}]}$ (resp. $w|_{(t_{cl}, T]}$) denotes the fault-on (resp. post-fault) trajectory of the selected dynamic/algebraic state~$w$.  

\begin{remark}
\textit{On observing the fault-on trajectory.} Observe that in~\eqref{eq:operator} the operator~$G^\dagger$ takes as the input the exact fault-on trajectory, \ie $w|_{[t_f,t_{cl}]}$. While this may be possible when using this operator setting for offline post-fault trajectory prediction (\ie using numerical simulation), in the online setting, it may not be possible to exactly capture~$w|_{[t_f,t_{cl}]}$ using PMUs. As a result, we propose next a relaxed operator regression problem for post-fault trajectory prediction.    
\end{remark}

We start by splitting the time domain of transient stability, denoted as~$[t_0,T]$, into two adjacent regions. The first region, denoted as~$T_u$, satisfies~$[t_f,t_{cl}] \subseteq T_u$ and corresponds to the time interval containing the domain of all possible fault-on trajectories. The second region, denoted as~$T_G$, satisfies $T_G = [0,T] \setminus T_u$ and corresponds to the time domain of post-fault trajectories. Thus, we define the ``relaxed'' operator~$G^\dagger$ for post-fault trajectory prediction as:
\begin{align} \label{eq:relaxed-operator}
    G^\dagger: w|_{T_u} \mapsto w|_{T_G}. 
\end{align}
To learn the operator~\eqref{eq:relaxed-operator}, we use the Deep Operator Network~(DeepONet) framework introduced in~\cite{lu2021learning} based on the universal approximation theorem to nonlinear operators~\cite{chen1995universal}.

\subsection{Review of DeepONet} \label{sub-sec:deeponet}
\begin{figure}[t!]
\centering
\includegraphics[width=.50\textwidth, height=5.5cm]{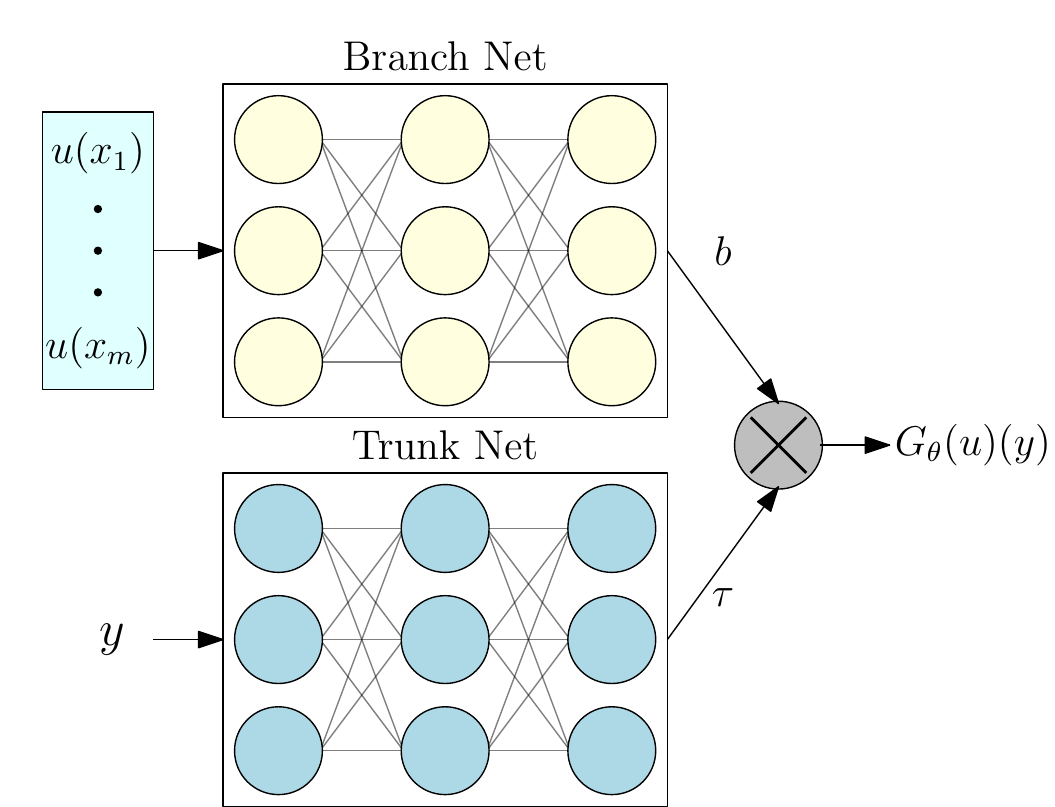}
\caption{The vanilla DeepONet architecture. The crossed node on the right indicates the DeepONet output, which we obtain by taking the inner product between the output features of the Branch~($b$) and Trunk~($\tau$) Nets.}
\label{fig:vanilla-deeponet}
\end{figure}

We now review the original formulation of DeepONet~\cite{lu2021learning}. Let $G^\dagger$ denote a nonlinear operator. This operator~$G^\dagger$ maps an input function~$u$ to an output function~$G^\dagger(u)$. Let $y \in Y$ denote a point in the output function domain $Y$ (often a subset of $\mathbb{R}^d$). Then, the goal of the DeepONet~$G_\theta$, with trainable parameters~$\theta \in \mathbb{R}^p$, is to approximate the operator~$G^\dagger(u)(y)$ at $y \in Y$. To this end, the DeepONet~$G_\theta$ uses the neural network architecture depicted in Figure~\ref{fig:vanilla-deeponet}, consisting of two sub-networks referred to as the Branch Net and the Trunk Net.      

The \textit{Branch}~Net processes the input function information. Let $(x_1,\ldots,x_m)$ denote points in the domain of~$u$ (we will refer to these points as the \textit{sensors}), such that $(u(x_1), \ldots, u(x_m))$ is a discrete representation of the input~$u$. The Branch Net takes this discretized~$u$ as the input and outputs a vector of features~$b \in \mathbb{R}^q$. On the other hand, the \textit{Trunk}~Net processes points in the domain~$Y$ of the output function. To this end, the Trunk Net takes~$y \in Y$ as the input and outputs a vector of features~$\tau \in \mathbb{R}^q$. Note that since the Trunk Net's output~$\tau$ solely depends on the input coordinates~$y$, it is natural to interpret the components of~$\tau$ as a collection of basis functions defined on~$Y$, \ie
$$
\tau = (\varphi_1(y), \ldots, \varphi_m(y)).
$$
The output of the DeepONet then combines the output features from the Branch Net~$b$ and the Trunk Net~$\tau$ using an inner product:
\begin{align} \label{eq:deeponet}
    G_\theta \left(u(x_1),\ldots,u(x_m)\right)(y) := \langle b,\tau \rangle + \tau_o= \sum_{i=1}^{q} b_i \cdot \varphi_i(y) + \tau_o,
\end{align}
where~$\tau_o$ is a trainable parameter added to the DeepONet's output. From the above, one can interpret the output of the Branch Net~$b$ as the trainable coefficients for the trainable basis functions~$\tau$ produced by the Trunk Net. To simplify our notation, in the rest of this work, we omit writing the DeepONet explicit dependency on the discretized input and use the simplified notation~$G_\theta(u)(y)$.

Finally, we train the DeepONet~$G_\theta$ to approximate the nonlinear solution operator~$G^\dagger$ by minimizing a mean square loss function on the training dataset $\mathcal{D} = \left\{u_i,y_i,G^\dagger(u_i)(y_i)\right\}_{i=1}^N$, \ie
\begin{align} \label{eq:deeponet-loss}
    \mathcal{L}(\theta) = \frac{1}{N} \sum_{i = 1}^N \left|G_\theta(u_i)(y_i) - G^\dagger(u_i)(y_i)\right|^2.
\end{align}

\begin{figure}[t!]
\centering
\includegraphics[width=.75\textwidth, height=6.0cm]{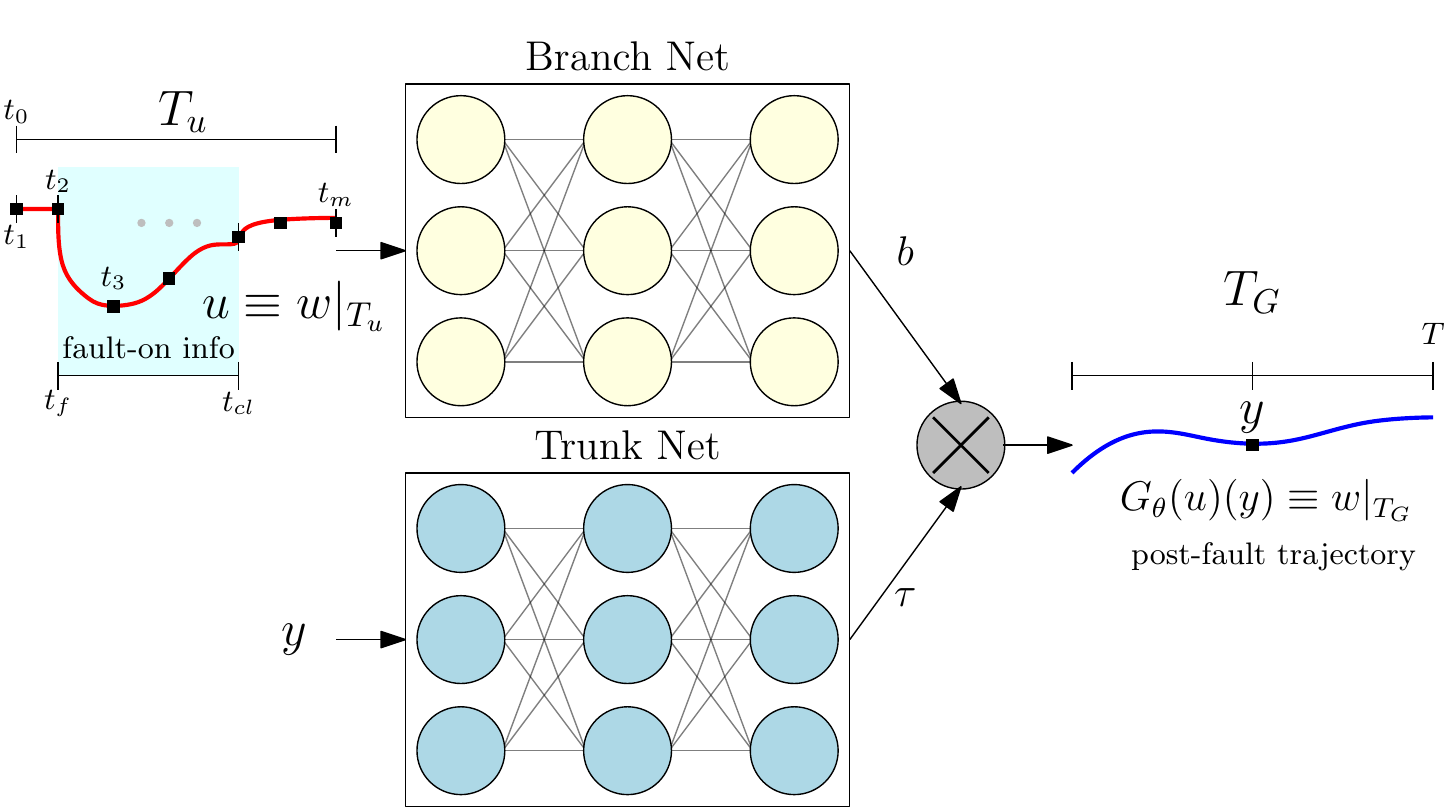}
\caption{The DeepONet for post-fault trajectory prediction. The sensors are $(t_1, \ldots, t_m)$, where $t_i \in T_u$. The input function is $u = w|_{T_u}$. The output function is $G^\dagger(u) = w|_{T_G}$. Thus the goal of DeepONet~$G_\theta$ is to approximate $G^\dagger(u)(y)$ at a point $y \in T_G$.}
\label{fig:deeponet-post-fault}
\end{figure}

\textit{The DeepONet for post-fault trajectory prediction:} Under the notation introduced in Section~\ref{sub-sec:deeponet}, the DeepONet that approximates the operator~\eqref{eq:relaxed-operator} for post-fault trajectory prediction (see Figure~\ref{fig:deeponet-post-fault}) is given as follows. The sensors $(x_1, \ldots, x_m)$ correspond to sampling times $(t_1, \ldots, t_m)$ within the time domain~$T_u$, \ie $t_i \in T_u$ for all $i=1,\ldots,m$. The input function corresponds to the trajectory containing fault-on information, \ie $u \equiv w|_{T_u}$. Moreover, the output function corresponds to the post-fault trajectory, \ie $G^\dagger(u) \equiv w|_{T_G}$. Thus, the domain of the output function corresponds to the time domain of post-fault trajectories, \ie $Y \equiv T_G$, and $y$ denotes a time sampled from~$T_G$. Throughout the rest of this paper, we keep the notation introduced in Section~\ref{sub-sec:deeponet} that follows~\cite{lu2021learning}. This is because the methods that we develop next for trustworthy post-fault trajectory prediction can also be used for other operator learning problems.  

We conclude this section with the following remark. DeepONets have shown remarkable approximation and generalization capabilities. Furthermore, DeepONets can learn operators efficiently with relatively small datasets. However, DeepONets inherits some of the problems from more traditional neural networks. Their reliability deteriorates when the input training data are noisy or do not cover all possible scenarios. One should note that noisy and incomplete training data is inherent in all power grid prediction problems. To overcome this reliability issue, we develop in the following sections two methods that enable DeepONets to provide trustworthy predictions through uncertainty quantification.

\section{Bayesian Deep Operator Network} \label{sec:bayesian-deeponet}
\begin{figure}[t!]
\centering
\includegraphics[width=.55\textwidth, height=7.5cm]{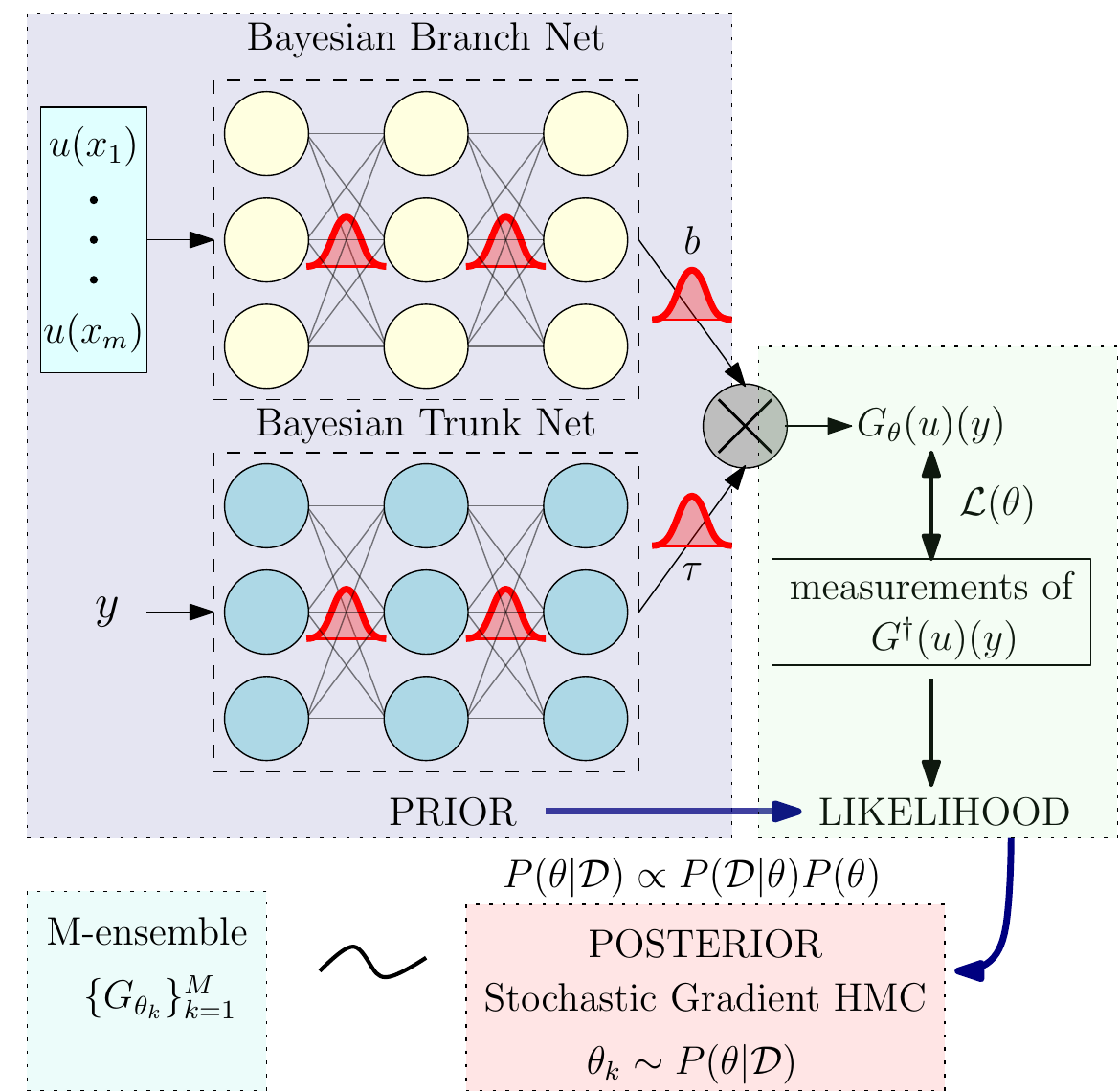}
\caption{The Bayesian DeepONet framework. The Bayesian DeepONet represents the prior for the trainable parameters~$\theta \in \mathbb{R}^p$. One computes the likelihood using the DeepONet forward pass and the training targets. Finally, one estimates the posterior using stochastic gradient Hamiltonian Monte-Carlo~(HMC).}
\label{fig:bayesian-deeponet}
\end{figure}

In this section, we propose a \textit{Bayesian DeepONet}~(B-DeepONet) framework (see Figure~\ref{fig:bayesian-deeponet}) to predict post-fault trajectories and quantify the aleatoric uncertainty. The proposed B-DeepONet consists of three components. (1) We represent the \textit{prior}~$p(\theta)$ of the trainable parameters~$\theta \in \mathbb{R}^p$ using the Bayesian DeepONet; (2) We compute the \textit{likelihood}~$p(\mathcal{D}|\theta)$ using the forward pass of DeepONets and the measurements $G^\dagger(u)(y)$; and (3) we sample from the \textit{posterior}~$p(\theta|\mathcal{D})$ using the (theoretically sound) Stochastic Gradient Hamiltonian Monte-Carlo algorithm. We review next the Hamiltonian Monte-Carlo~(HMC) algorithm~\cite{neal2011mcmc}.   
\subsection{Hamiltonian Monte-Carlo~(HMC)} \label{sub-sec:HMC}
The goal of HMC~\cite{neal2011mcmc} is to sample from the posterior distribution of the trainable parameters~$\theta$ given a set of independent observations~$d \in \mathcal{D}$. This target posterior is:
$$
p(\theta | \mathcal{D}) \propto \exp(-U(\theta)),
$$
where $U$ is the potential energy defined as
\begin{align} \label{eq:potential}
    U(\theta) = \underbrace{- \sum_{d \in \mathcal{D}} \log p(d | \theta)}_{\text{likelihood}} - \underbrace{\log p(\theta)}_{\text{prior}}.
\end{align}

To generate the samples, HMC uses the Hamiltonian system:
$$
H(\theta, r) = U(\theta) + \frac{1}{2} r^\top M^{-1} r.
$$
Here, $r \in \mathbb{R}^p$ is a vector of auxiliary momentum variables, and M (often set to be the identity matrix) is the mass matrix. The samples generated using~$H(\theta, r)$ follow the joint distribution of~$(\theta,r)$ given by
$$
\pi(\theta, r) \propto \exp \left(- H(\theta,r) \right).
$$
If we discard the obtained $r$ samples, then the $\theta$ samples have as marginal distribution the desired target distribution~$p(\theta | \mathcal{D})$~\cite{neal2011mcmc}. To propose new samples, HMC simulates the following Hamiltonian dynamics:
\begin{subequations}\label{eq:hamiltonian-dynamics}
\begin{align} 
    d\theta &= M^{-1} r~dt \\
    dr &= - \nabla U(\theta)~dt.
\end{align}
\end{subequations}
To simulate from ~\eqref{eq:hamiltonian-dynamics}, the HMC algorithm uses, for \eg the \textit{leapfrog} discretization method. In addition, a \textit{Metropolis-Hasting} step is used to accept/reject samples, which can contain inaccuracies as a result of the leapfrog discretization.

The main disadvantage of HMC is the computational cost we require to (1) compute the gradient~$\nabla U(\theta)$ when the dataset~$\mathcal{D}$ is large and (2) perform the Metropolis-Hasting step. To tackle this computational cost drawback, we propose using \textit{Stochastic Gradient} HMC~\cite{chen2014stochastic} for {\color{black}sampling} the posterior distribution of the DeepONet parameters~$\theta$.

\subsection{Stochastic Gradient HMC~(SGHMC)} \label{sub-sec:SGHMC}
In SGHMC~\cite{chen2014stochastic}, instead of computing the expensive gradient~$\nabla U(\theta)$, which requires traversing the entire dataset~$\mathcal{D}$, SGHMC considers a noisy estimate of this gradient based on a mini-batch $\tilde{\mathcal{D}} \subset \mathcal{D}$ sampled uniformly at random,  \ie
\begin{align} \label{eq:noisy-gradient}
\nabla \tilde{U} (\theta) = - \frac{|\mathcal{D}|}{|\tilde{\mathcal{D}}|} \sum_{d \in \tilde{\mathcal{D}}} \nabla \log (d | \theta) - \nabla \log p(\theta).
\end{align}
If the $d$'s are independent, then via the central limit theorem, the noisy gradient is given approximately by~\cite{chen2014stochastic}:
$$
\nabla \tilde{U}(\theta) \approx \nabla U(\theta) + \mathcal{N}(0, V(\theta)),
$$
where $V(\theta)$ is the covariance of the stochastic gradient noise,\footnote{We have adopted the same notation as in~\cite{chen2014stochastic}, that is $\mathcal{N}(\mu, \Sigma)$ represents a random variable distributed according to the multivariate Gaussian.} which could depend on~$\theta$ and $|\tilde{\mathcal{D}}|$.

To simulate samples using SGHMC, we use the following modified version of the dynamics~\eqref{eq:hamiltonian-dynamics} that adds a \textit{friction} term to the momentum update, \ie
\begin{subequations}\label{eq:hamiltonian-dyn-with-friction}
\begin{align} 
    d\theta &= M^{-1} r~dt \\
    dr &= - \nabla \tilde{U}(\theta)~dt-CM^{-1}rdt + \mathcal{N}(0,2 (C-\hat{B})dt) + \mathcal{N}(0, 2Bdt).
\end{align}
\end{subequations}
Here $C \succeq \hat{B}$ is a user-specified \textit{friction} term and $\hat{B}$ is an estimate of the noise model $B(\theta) = \frac{1}{2 \epsilon} V(\theta)$, where $\epsilon$ is the step size of the discretization scheme. As demonstrated in~\cite{chen2014stochastic}, the above dynamics, with the proper selection of $C$ and $\hat{B}$, yields the desired stationary distribution $\pi(\theta, r) \propto e^{-H(\theta, r)}$. We provide the details for implementing SGHMC in Algorithm~\ref{alg:SGHMC}.
\begin{algorithm}[t]
\DontPrintSemicolon
\SetAlgoLined
\textbf{Require:} initial state~$\theta^0$ and time step size~$\epsilon_t$;
\For{$k = 1,\ldots,N$}{
  sample $r^{k-1} \sim \mathcal{N}(0,M)$\;
  $(\theta_0,r_0) = (\theta^{k-1}, r^{k-1})$\;
  \For{$i=1,\ldots,m$}{
  simulate equation~\eqref{eq:hamiltonian-dyn-with-friction}:\;
  $\theta_i = \theta_{i-1} + \epsilon_t M^{-1} r_{i-1}$\;
  $r_i = r_{i-1} - \epsilon_t \nabla \tilde{U}(\theta_i) - \epsilon_t C M^{-1} r_{i-1} + \mathcal{N}(0, 2(C-\hat{B})\epsilon_t)$\;
  }
  $(\theta^k,r^k) = (\theta_m, r_m)$\;
  }
  Calculate $\{\{G_{\theta^{N+1-k}}(u)(y) : y \in Y_\text{mesh} \}_{k=1}^{M}$ as predicted samples of the true post-fault trajectory $\{G^\dagger(u)(y) : y \in Y_\text{mesh}\}$\;
 \caption{Stochastic Gradient HMC}
 \label{alg:SGHMC}
\end{algorithm}

\begin{figure}[t!]
\centering
\includegraphics[width=.60\textwidth, height=7.5cm]{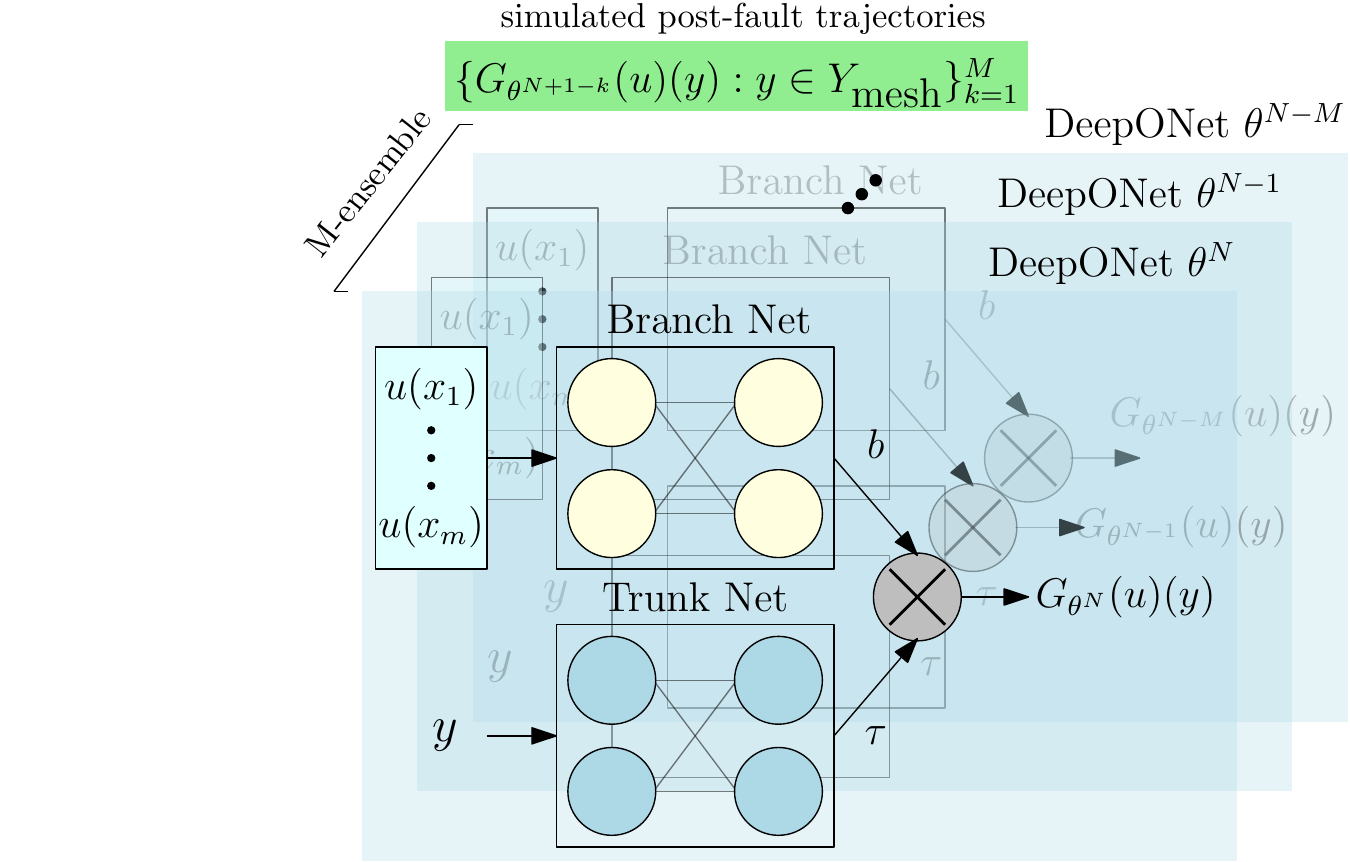}
\caption{$M$-ensemble of DeepONet parameters~$\{\theta^{N+1 - k}\}$ sampled from the posterior~$p(\theta|\mathcal{D})$ using SGHMC. We use this ensemble to simulate post-fault trajectories $\{G_{\theta^{N+1-k}}(u)(y):y \in Y_\text{mesh}\}_{k=1}^M$. The mean of these trajectories is used to predict the true trajectory $\{G^\dagger(u)(y): y \in Y_\text{mesh}\}$. The {\color{black}standard} deviation is used to quantify the aleatoric uncertainty.}
\label{fig:ensemble}
\end{figure}

To provide a trustworthy prediction of the post-fault trajectories (see Algorithm~\ref{alg:SGHMC} for details), denoted as $\{G^\dagger(u)(y) : y \in Y_\text{mesh}\}$, for a given input~$u$ and over the given mesh of points $y \in Y_\text{mesh} \subset Y \equiv T_G$, we must sample from the posterior~$p(\theta|\mathcal{D})$ and obtain the $M$-ensemble of DeepONets with parameters~$\{\theta^{N+1-k}\}_{k=1}^M$ depicted in Figure~\ref{fig:ensemble}. Then, using this ensemble, we can simulated $M$ trajectories~$\{G_{\theta^{N+1-k}}(u)(y):y \in Y_\text{mesh} \}_{k=1}^M$ and compute their statistics. Here, we compute the mean and standard deviation. We use the mean to predict the true solution $\{G^\dagger(u)(y) : y \in Y_\text{mesh}\}$ and the standard deviation to quantify the uncertainty. 

One should note that trustworthy post-fault trajectory prediction requires simulating $M$ trajectories $\{G_{\theta^{N+1-k}}(u)(y) : y \in Y_\text{mesh}\}_{k=1}^M$ using the $M$-ensemble of DeepONets produced by Algorithm~\ref{alg:SGHMC}. However, simulating these many trajectories for large-scale power grids might become computationally expensive,  which precludes the application of the proposed B-DeepONet to online post-fault trajectory prediction. While some of the computational costs may be alleviated using a parallel computing framework, we believe B-DeepONets are more suitable for offline post-fault trajectory prediction. 

\section{Probabilistic Deep Operator Network} \label{sec:probabilistic-deeponet}
The main drawback of most data-driven methods for transient stability is that they fail to handle the inherent trade-off between providing (1) trustworthy predictions and (2) efficient/fast predictions. Our proposed B-DeepONet (see Section~\ref{sec:bayesian-deeponet}) gave us the tools (supported by the solid mathematical results detailed in~\cite{chen2014stochastic}) to provide reliable predictions by estimating their aleatoric uncertainty. However, the increased predictive reliability gained by B-DeepONets comes at the price of decreased predictive efficiency. 

Motivated by the need for both trustworthy and efficient post-fault trajectory predictions, here, in this section, we take a more \textit{empirical} approach and design a probabilistic training framework for DeepONets, which we refer to as \textit{Prob-DeepONet}. Prob-DeepONet was introduced in one of the authors' previous works~\cite{winovich2019convpde} and will enable us to provide efficient point-wise aleatoric uncertainty estimates of post-fault trajectory predictions.

\begin{figure}[t!]
\centering
\includegraphics[width=.65\textwidth, height=6.0cm]{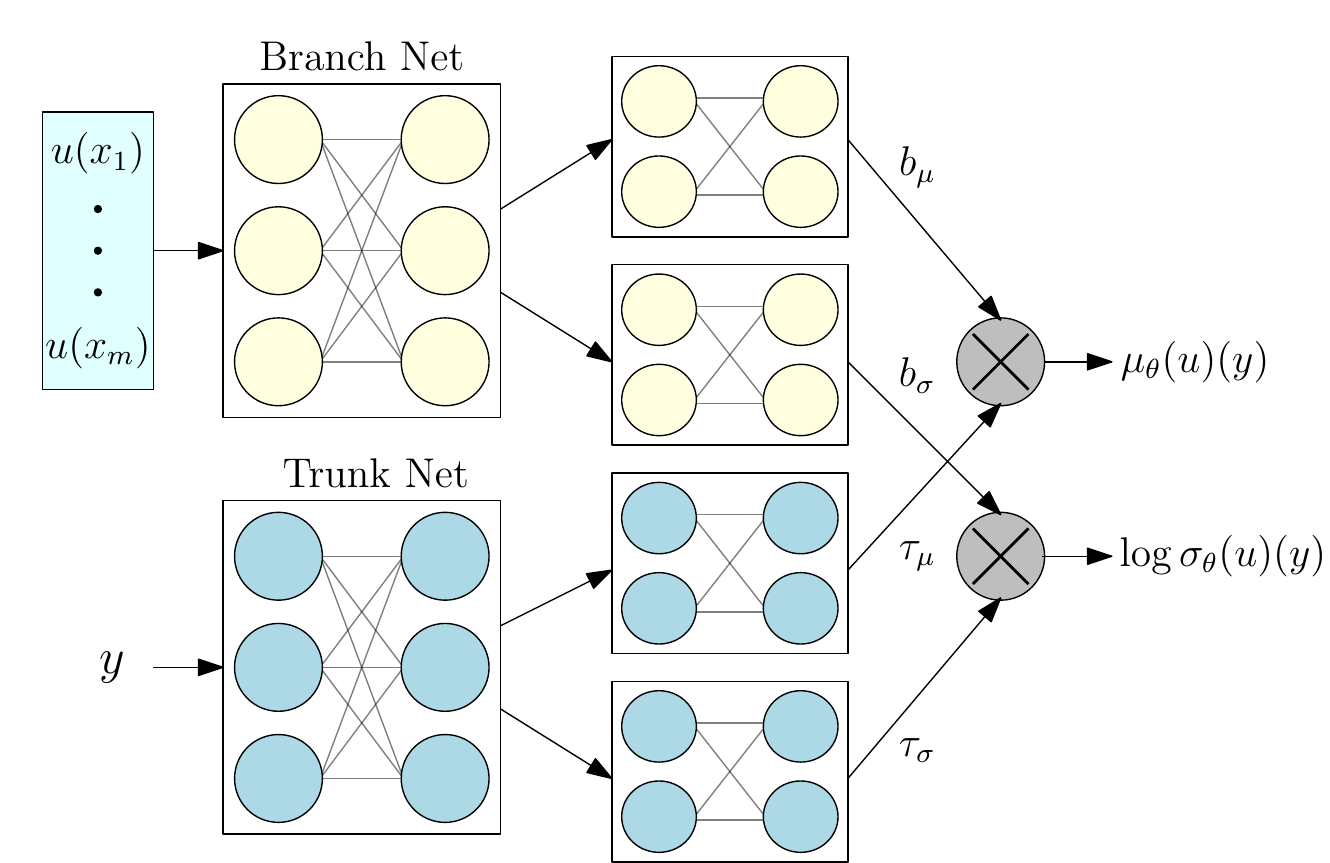}
\caption{The Probabilistic DeepONet framework. The DeepONet outputs are interpreted as parameters of a Gaussian distribution~$\mathcal{N}(\mu_\theta(u)(y), \sigma^2_\theta(u)(y))$ of the \textit{true} values $G^\dagger(u)(y)$. The loss function is given by the negative log-likelihood~\eqref{eq:NLL}.}
\label{fig:probabilistic-deeponet}
\end{figure}

Prob-DeepONet requires building a deep operator network architecture (depicted in Figure~\ref{fig:probabilistic-deeponet}) with the ability to learn estimates of aleatoric uncertainty. As described in~\cite{winovich2019convpde}, this can be achieved by casting the DeepONet outputs as Gaussian posterior distributions. We select the outputs to represent Gaussian distributions because we have observed (see Appendix~1 for more details) that the \text{network} errors ($G_\theta(u)(y) - G^\dagger(u)(y)$) experienced by DeepONets when trained via gradient descent algorithms (\eg Adam) follow a Gaussian distribution.  More specifically, we build a DeepONet architecture that predicts the point-wise statistics~$\mu_\theta(u)(y)$ and $\sigma_\theta(u)(y)$ of the Gaussian distribution $\mathcal{N}(\mu_\theta(u)(y), \sigma^2_\theta(u)(y))$ of the \textit{true} target operator values~{\color{black}$G^\dagger(u)(y)$}.

\textit{Network architecture.} To accommodate the probabilistic prediction for DeepONets, we split the final layers of both the Branch Net and the Trunk as depicted in Figure~\ref{fig:probabilistic-deeponet}. These new layers provide independent processing for the mean~$\mu_\theta$ and the standard deviation~$\sigma_\theta$ of the DeepONet's Gaussian estimates of the output. To avoid numerical instability, we follow~\cite{winovich2019convpde} and construct the DeepONet to predict the $\log$ standard deviation $\log \sigma_\theta(u)(y)$ and recover the standard deviation using the exponential function, \ie $\sigma_\theta(u)(y) = \exp(\log \sigma_\theta(u)(y))$.

\textit{Probabilistic Training.} During training, our goal is to calibrate the DeepONet's parameters~$\theta$ to produce the Gaussian estimates of predictive uncertainty. This requires designing a probabilistic training protocol to maximize the associated likelihood of observing the target post-fault trajectory values. Equivalently, one can minimize the negative-log likelihood of the target post fault trajectories, \ie  
\begin{align} \label{eq:NLL}
    \mathcal{L}_p(\theta) = \frac{1}{N} \sum_{i = 1}^N \frac{1}{2} \frac{(\mu_\theta(u_i)(y_i) - G^\dagger(u_i)(y_i))^2}{ \sigma^2_\theta(u_i)(y_i)} + \frac{1}{2} \log(2 \pi \sigma^2_\theta(u_i)(y_i)).
\end{align}
As described in~\cite{winovich2019convpde}, during probabilistic training, this negative-log likelihood loss guides the DeepONet to build trustworthiness gradually by decreasing the predicted standard deviation~$\sigma_\theta$. 

\textit{Trustworthy prediction of post-fault trajectories.} To predict and estimate the aleatoric uncertainty of a post fault trajectory~$\{G^\dagger(u)(y):y \in Y_\text{mesh}\}$, Prob-DeepONet predicts the mean trajectory $\{\mu_\theta(u)(y): y \in Y_\text{mesh}\}$ and the standard deviation $\{\sigma_\theta(u)(y): y \in Y_\text{mesh}\}$. The mean trajectory represents the prediction of the true trajectory $\{G^\dagger(u)(y):y \in Y_\text{mesh}\}$ and we use the standard deviation to construct an estimate of its uncertainty. Let us conclude this section by noting that compared to B-DeepONet, which requires $M$ DeepONet forward passes to construct an estimate of the uncertainty (see Figure~\ref{fig:ensemble}), Prob-DeepONet estimates uncertainty instantaneously at virtually no extra cost, \ie it uses only \textit{one} Prob-DeepONet forward pass. In the next section, we use a series of numerical experiments to demonstrate the B-DeepONet and Prob-DeepONet's effectiveness in providing trustworthy predictions.   

\section{Numerical Experiments} \label{sec:numerical-experiments}
This section presents four numerical experiments aimed at demonstrating the predictive power and reliability of the proposed B-DeepONet and Prob-DeepONet. The first experiment~(Section~\ref{sub-sec:predictive-power}) compares the predictive power of B-DeepONet, Prob-DeepONet, and vanilla DeepONet. The second experiment~(Section~\ref{sub-sec:reliability}) builds confidence intervals to verify how reliable the predictions of B-DeepONet and Prob-DeepONet are. The third experiment~(Section~\ref{sub-sec:robustness}) tests the performance of the proposed methods when the input fault-on trajectory is contaminated with noise. The final experiment~(Section~\ref{sub-sec:alarms}) tests how well the proposed methods predict violations of pre-defined voltage limits. We start this section by describing the power grid simulation framework used to generate the datasets.

\subsection{Dataset generation}  \label{sub-sec:dataset-generation}
\begin{figure}[t!]
\centering
\includegraphics[width=.70\textwidth, height=6.5cm]{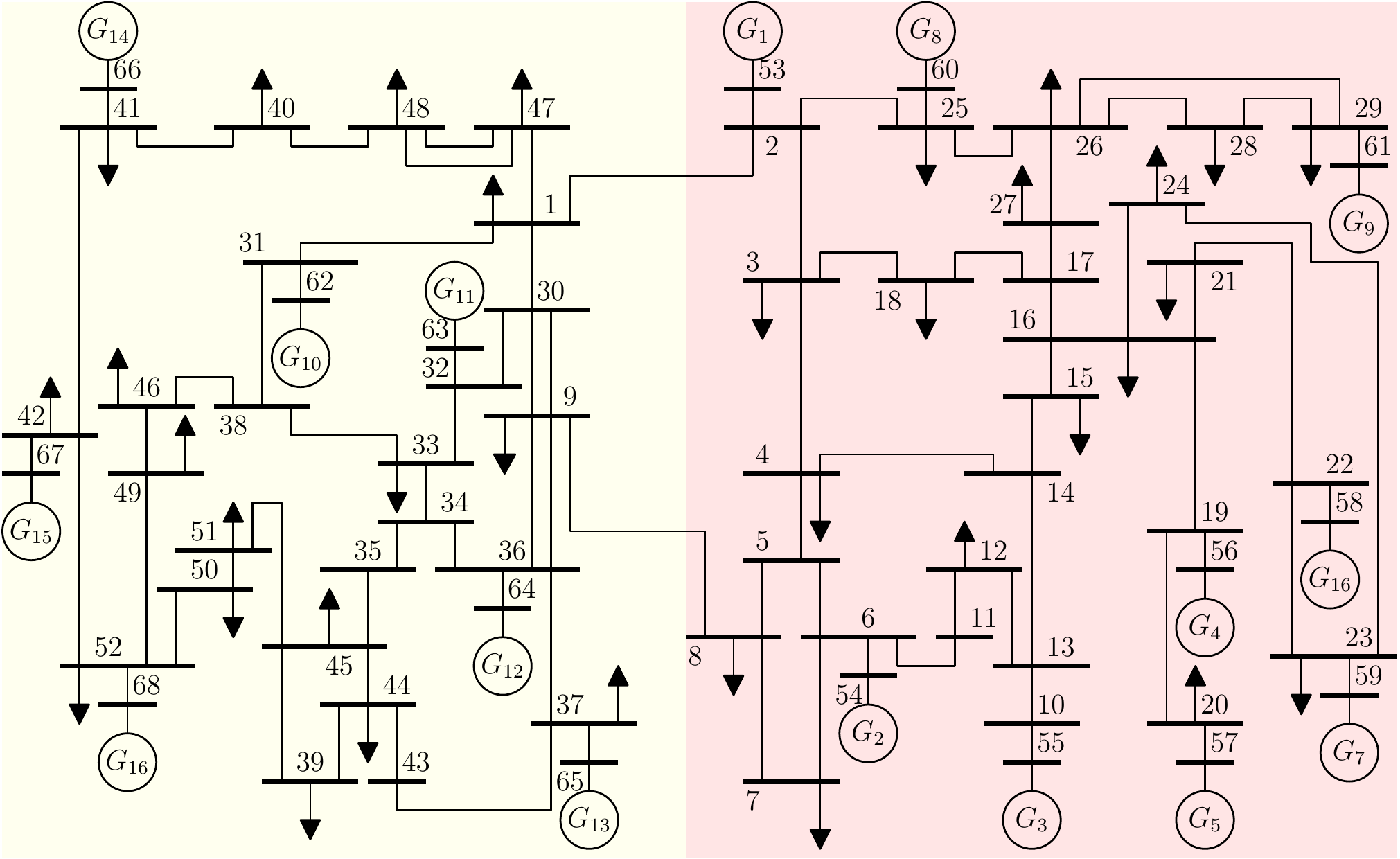}
\caption{One-line diagram of the New York-New England 16-generator 68-bus power grid.}
\label{fig:68-bus-system}
\end{figure}

We used time-domain simulations on the New York-New England 16-generator 68-bus power grid model (see Figure~\ref{fig:68-bus-system}) to generate the training and testing datasets. We performed the simulations using EPTOOL --a tool developed based on the Power System Toolbox~\cite{chow1992toolbox}. We also considered the full dynamical model of the generators, including turbine governors, excitation systems, and power system stabilizer.

In the numerical experiments, we learn the operator~$G^\dagger$ that predicts the post-fault trajectories of the voltage magnitude of bus~$h=19$ (see Figure~\ref{fig:68-bus-system}). However, we remark that the proposed framework can be implemented on other buses or states, \eg frequency.  We leave for our future work the design of a transfer learning strategy to reduce the training time of DeepONets for other buses and states (\eg the frequency at some bus~$j$) using a trained DeepONet of a selected state (\eg the voltage of bus~$h$, where $h \neq j$).

To train and test the B-DeepONet and Prob-DeepONet, we generate two data pools of transient response trajectories. We generate the first data pool of $1000$ trajectories by simulating the grid's transient response to $N-$1 faults. Similarly, we generate the second data pool of $1000$ trajectories by simulating $N-$2 faults. Each trajectory is simulated for $T=9$ seconds with a resolution of $100$ [Hz]. As a result, all the $N= 2000$ trajectories have $N_\text{samples} = 900$ samples.

To simulate faults, we proceed as follows. We assume the power grid operates on a stable equilibrium before the fault, \ie within the interval $[0,t_f)$. At the time $t_f$, the grid experiences the fault. We simulate $N-$1 (resp. $N-$2) faults by disconnecting uniformly at random one (resp. two) transmission lines. Then, at time $t_{cl}$, we assume the reclosing system clears the fault. Finally, the post-fault trajectory occurs within the interval $(t_{cl},T]$. 

To simulate each fault, we must define the values for~$t_f$ and~$t_{cl}$. Here, for all $N$ trajectories, we fix the value of the clearing time to $t_{cl} = 2$ seconds. Let~$\Delta t_f$ denote the duration time of the fault, \ie $\Delta t_f:= t_{cl} - t_f$. We simulate $N$ $\Delta t_f$'s from the uniform distribution $\mathcal{U}(0.2,0.5)$. Then, for each simulated~$\Delta t_f$, we obtain $t_f$ via $t_f = 2.0 - \Delta t_f$.

Finally, to generate the training and testing data pools, we combine and shuffle the $N-$1 and $N-$2 data pools and perform a $70-30 \%$ split.

\textit{Training dataset~$\mathcal{D}_\text{train}$.} Here, we use the training data pool of $N_\text{train}$ transient voltage trajectories to construct our training dataset~$\mathcal{D}_\text{train}$. We first split the transient simulation domain~$[0, T]$ into the two adjacent domains. $T_u = [0, t_{cl}]$ and $T_G=(t_{cl}, T]$. $T_u$ contains the fault-on information and $T_G$ represents the domian of post-fault trajectories.

The Branch Net~$u$ takes as inputs the voltage trajectories in the interval~$T_u$. We discretize these inputs using $m=200$ sensors. Thus, the discretized input $(u(x_1),\ldots,u(x_m))$ coincides with the voltage trajectory simulated at a resolution of 100~[Hz]. One can also use a different number of sensors. However, this requires performing interpolation.  

On the other hand, the Trunk Net takes as inputs time values~$y$ sampled uniformly from the interval~$T_G$. Finally, the target operator values~$G^\dagger(u)(y)$ for supervised training correspond to the transient voltage response evaluated at time~$y$ for the given input (fault-on trajectory)~$u$. Note that if $y$ does not coincide with the simulation resolution of~$100$ [Hz], we need to perform interpolation. The training dataset is then
$$
\mathcal{D}_\text{train} = \{(u_i(x_1),\ldots,u_i(x_m)), y_i, G^\dagger(u_i)(y_i)\}_{i=1}^{N_\text{train}}.
$$
In the above~$\mathcal{D}$, we have shown the case when  only $Q=1$ point~$y_i \in T_G$ is sampled per input trajectory~$u_i$. 
However, we remark that we can always sample multiple points~$\{y_i^j\}_{j=1}^Q$ with $Q > 1$ per input trajectory~$u_i$. In all our experiments, we sampled~$Q=10$ points.

\textit{Testing dataset~$\mathcal{D}_\text{test}$.} The testing dataset consists of all the $N_\text{test}$ simulated trajectories of the testing data pool. Using the DeepONet notation, the post-fault component of a testing trajectory is denoted as $\{G(u)(y):y \in Y_\text{mesh}\}$. This post-fault component corresponds to the transient trajectory evaluated at points~$y \in Y_\text{mesh} \subset T_G$ for the given input~$u$. In all our experiments, we selected $Y_\text{mesh}$ to be the set of 500 linearly spaced points within the post-fault interval $T_G$.  Finally, to evaluate the performance of B-DeepONet and Prob-DeepONet on the testing dataset, we use two metrics: the $L_1$ and the $L_2$ relative errors. 

\subsection{Training protocols and neural networks}
\textit{Training protocols.} We implemented Prob-DeepONet and B-DeepONet using PyTorch and published the code in GitHub.\footnote{The code will be available after publication.} To train the Prob-DeepONet and vanilla DeepONet, we used the Adam~\cite{kingma2014adam} optimizer with default hyper-parameters and initial learning rate $\eta = 10^{-4}$. We reduced this learning rate whenever the loss reached a plateau or increased. For B-DeepONet, we manually calibrated all the parameters of SGHMC (Algorithm~\ref{alg:SGHMC}). 

\textit{Neural networks.} We implemented the Branch and the Trunk Nets using the modified fully-connected architecture proposed in ~\cite{wang2020understanding} and used in our previous paper~\cite{moya2021dae}. The forward pass of this modified network is 
\begin{align*}
&U = \phi(X W^1 + b^1),~V = \phi(X W^2 + b^2) \\
&H^{(1)} = \phi(X W^{z,1} + b^{z,1}) \\
&Z^{(k)} = \phi(H^{(k)} W^{z,k} + b^{z,k}),~k=1, \ldots,d \\
&H^{(k+1)} = (1 - Z^{(k)}) \odot U + Z^{(k)} \odot V,~k=1,\ldots,d \\
&f_\theta (x) = H^{(d+1)} W +b,
\end{align*}
Here, $X$ is the input tensor to the neural network, $d$ is the number of hidden layers (\ie the network's depth), $\odot$ is the element-wise product, and $\phi$, in this paper, is a point-wise sinusoidal activation function. The trainable parameters are
$$
\theta = \{W^1, b^1, W^2, b^2, \{W^{z,l}, b^{z,l}\}_{l=1}^d, W, b\}.
$$
We initialized all these parameters using the Glorot uniform algorithm.

\subsection{Experiment 1: Mean prediction and generalization} \label{sub-sec:predictive-power}
This experiment compares the mean predictive power and the generalization capability of B-DeepONet, Prob-DeepONet, and vanilla DeepONet. To this end, we train the models for 10000 epochs. Once we complete training, we select the best-trained models for vanilla DeepONet and Prob-DeepONet, and the $M-$ensemble produced by B-DeepONet. We use these trained models to compute the average and standard deviation~(st.dev.) of the $L_1-$ and $L_2-$ relative errors for 100 test trajectories selected uniformly at random from~$\mathcal{D}_\text{test}$.

We report the results in Table~\ref{table:predictive-power}. These results show that B-DeepONet and Prob-DeepONet provide similar mean predictive power and generalization capability to the vanilla DeepONet. We believe this is a remarkable result because our proposed methods were not designed to provide improved mean prediction of post-fault trajectories. 

\begin{table}[ht]
\centering
\begin{tabular}{| c | c  c  c|}
\hline
 & B-DeepONet & Prob-DeepONet & vanilla DeepONet\\
\hline
\hline
mean $L_1$ & 1.81  \% & 1.62  \% & 1.63  \% \\
st.dev. $L_1$ & 1.88 \% & 2.02 \% & 2.0 \%\\
\hline
mean $L_2$ & 2.23 \% & 1.95 \% & 2.08 \%\\
st.dev. $L_2$ & 2.40 \% & 2.51 \% & 2.57 \% \\
\hline
\end{tabular}
\caption{The average and standard deviation (st.dev.) of the $L_1-$ and $L_2-$ relative errors of 100 trajectories selected uniformly at random from $\mathcal{D}_\text{test}$.}
\label{table:predictive-power}
\end{table}

\subsection{Experiment 2:~Reliability} \label{sub-sec:reliability}
This experiment verifies how well the M-ensemble produced by B-DeepONet and the trained Prob-DeepONet quantify the uncertainty. To this end, we construct a $95\%$ confidence interval~(CI) for post-fault trajectory prediction. B-DeepONet constructs the confidence interval by computing the standard deviation over the $M$ simulated post-fault trajectories $\{G_{\theta^{N+1-k}}(u)(y): y \in Y_\text{mesh}\}_{k=1}^M$. Prob-DeepONet constructs the confidence interval using the point-wise predicted standard deviation of the post fault trajectory, that is, $\{\sigma_{\theta^*}(u)(y):y \in Y_\text{mesh}\}$. 

To measure how reliable is the predicted uncertainty, we define the following ratio~$\epsilon_\text{ratio}$:
$$
\epsilon_\text{ratio} = \frac{\# \text{ of points of the true post-fault trajectory}~\{G^\dagger(u)(y):y \in Y_\text{mesh}\} \text{ in the CI}}{\# \text{ of points of the true post-fault trajectory}~\{G^\dagger(u)(y):y \in Y_\text{mesh}\}} \cdot 100 \%
$$
Note that $\epsilon_\text{ratio}$ measures the percentage of the true post-fault test trajectory captured by the predicted confidence interval.

\begin{figure}[t!]
\centering
\begin{subfigure}[b]{0.49\textwidth}
\centering
\includegraphics[width=1.0\textwidth, height=5.5cm]{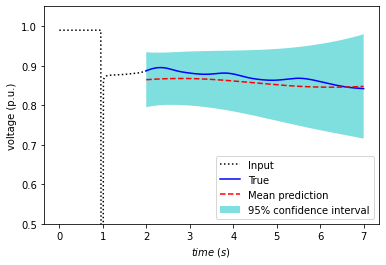}
\end{subfigure}
\begin{subfigure}[b]{0.49\textwidth}
\centering
\includegraphics[width=1.0\textwidth, height=5.5cm]{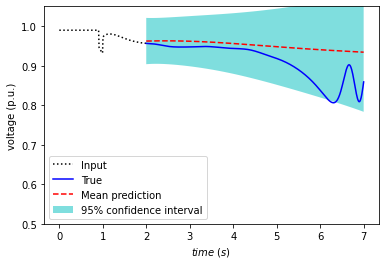}
\end{subfigure}
\\
\begin{subfigure}[b]{0.49\textwidth}
\centering
\includegraphics[width=1.0\textwidth, height=5.5cm]{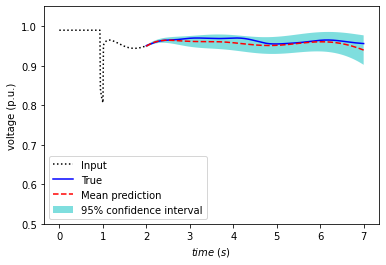}
\end{subfigure}
\begin{subfigure}[b]{0.49\textwidth}
\centering
\includegraphics[width=1.0\textwidth, height=5.5cm]{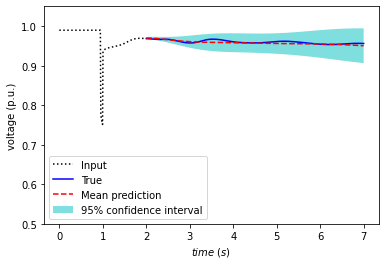}
\end{subfigure}
\caption{\textit{Top:} B-DeepONet mean prediction and confidence interval~(CI) for two post-fault trajectories selected at random from~$\mathcal{D}_\text{test}$. The CI contains $\epsilon_\text{ratio} = 100 \%$ of the true test trajectories. \textit{Bottom:} Prob-DeepONet mean prediction and confidence interval for two post-fault trajectories selected at random from~$\mathcal{D}_\text{test}$. The CI contains $\epsilon_\text{ratio} = 100 \%$ of the true test trajectories.}
\label{fig:confidence-intervals}
\end{figure}

Figure~\ref{fig:confidence-intervals} depicts the mean prediction and confidence intervals generated by the proposed methods for four test trajectories selected uniformly at random from~$\mathcal{D}_\text{test}$. One should note that the constructed CIs capture pretty well the true test post-fault trajectories, \ie their ratio is $\epsilon_\text{ratio} = 100 \%$. 

We also report in Table~\ref{table:ratios} the average ratio~$\epsilon_\text{ratio}$ for 100 trajectories selected from~$\mathcal{D}_\text{test}$. As expected, both B-DeepONet and Prob-DeepONet capture with their confidence intervals $\approx 95 \%$ of all the test trajectories. 

\begin{table}[h!]
\centering
\begin{tabular}{| c | c  c |}
\hline
 & B-DeepONet & Prob-DeepONet \\
\hline
\hline
$\epsilon_\text{ratio}$ & 93.32 \% & 94.27 \% \\
\hline
\end{tabular}
\caption{The average $\epsilon_\text{ratio}$ of all the post-fault trajectories in~$\mathcal{D}_\text{test}$.}
\label{table:ratios}
\end{table}

To conclude this experiment, we compute, for B-DeepONet and Prob-DeepONet, the fraction of points of the test trajectories falling within a specified standard deviation range. To this end, let us first define for B-DeepONets the predicted mean of the test trajectory~$\{G^\dagger(u)(y):y \in Y_\text{mesh}\}$ as
$\{\hat{\mu}_\theta^B(u)(y):y \in Y_\text{mesh}\}$ and the predicted standard deviation as $\{\hat{\sigma}^B_\theta(u)(y):y \in Y_\text{mesh}\}$. 

Formally, we aim to compute the average number of points where the true post-fault trajectory differs from the predictive mean by less than a factor~$\chi$ of the predicted standard deviation:
\begin{align}  \label{eq:average}
\frac{1}{|\mathcal{D}_{\text{test}}|} \frac{1}{|Y_\text{mesh}|} \sum_{u \in \mathcal{D}_\text{test}} \sum_{y \in Y_\text{mesh}} \mathbbm{1} \{E(u)(y)\},
\end{align}
where $\mathbbm{1}\{\cdot\}$ is the indicator function and $E(u)(y)$ is the aforementioned event:
$$
E(u)(y):=|\hat{\mu}(u)(y) - G^\dagger(u)(y)| \le \chi \cdot \hat{\sigma}(u)(y).
$$
Here, $\hat{\mu} \equiv \hat{\mu}_\theta^B$ and $\hat{\sigma} \equiv \hat{\sigma}_\theta^B$ for B-DeepONets, and $\hat{\mu} \equiv \mu_\theta$ and $\hat{\sigma} \equiv \sigma_\theta$ for Prob-DeepONet.

\begin{figure}[t!]
\centering
\includegraphics[width=0.5\textwidth, height=5.5cm]{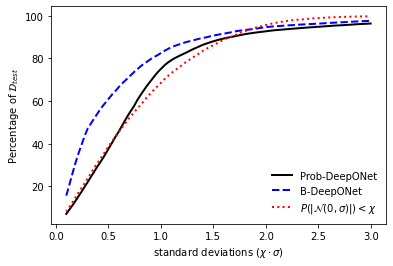}
\caption{Quantitative analysis of the B-DeepONet (\textit{left}) and Prob-DeepONet (\textit{right}) predicted standard deviation values for all post-fault trajectories in $\mathcal{D}_\text{test}$. We show the percentage of points where the difference between the mean prediction and true solution is less than a factor~$\chi$ of the predicted standard deviation. They both closely follow the cumulative distribution of the half-normal distribution.}
\label{fig:xi-analysis}
\end{figure}

Figure~\ref{fig:xi-analysis} compares the average points obtained
using~\eqref{eq:average} for B-DeepONet and Prob-DeepONet against  $\mathbb{P}(|\mathcal{N}(0,1)| \le \chi)$, \ie the average number of points of a Normal distribution lying within a factor~$\chi$ of the unit standard deviation. As expected, the Prob-DeepONet result remains close to the Normal Distribution. Such a result confirms our assumption of Gaussian distributed network errors. Remarkably, the B-DeepONet average remains also close to the Normal distribution average. This experimental evidence suggests that the proposed methods adequately measure the predictive uncertainty.

\subsection{Experiment 3: Robustness to noise} \label{sub-sec:robustness}
In this experiment, we verify the performance of the proposed methods when the test input~$u$, containing fault-on information, is contaminated with noise. This scenario can happen in the online setting, which uses noisy measurements collected using PMUs. 

We assume the noisy discretized inputs~$\tilde{u}(x_i)$, for $i=1,\dots,m$ are independently Gaussian distributed centered at the true input value, \ie
$$
\tilde{u}(x_i) = u(x_i) + \epsilon_i, \qquad i=1,\ldots,m,
$$
where $\epsilon_i$ is an independent Gaussian noise with zero mean and standard deviation~$\sigma$, \ie $\epsilon_i \sim \mathcal{N}(0, \sigma^2)$. We select the noise level to be~$\sigma = 0.01$.

\begin{figure}[t!]
\centering
\begin{subfigure}[b]{0.49\textwidth}
\centering
\includegraphics[width=1.0\textwidth, height=5.5cm]{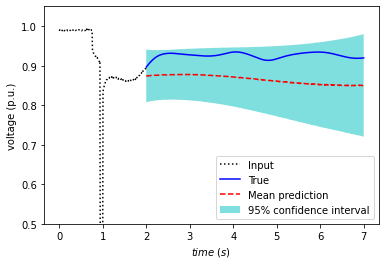}
\end{subfigure}
\begin{subfigure}[b]{0.49\textwidth}
\centering
\includegraphics[width=1.0\textwidth, height=5.5cm]{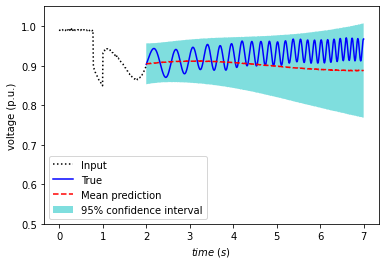}
\end{subfigure}
\\
\begin{subfigure}[b]{0.49\textwidth}
\centering
\includegraphics[width=1.0\textwidth, height=5.5cm]{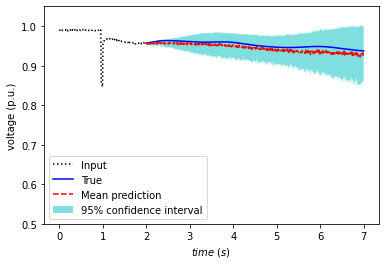}
\end{subfigure}
\begin{subfigure}[b]{0.49\textwidth}
\centering
\includegraphics[width=1.0\textwidth, height=5.5cm]{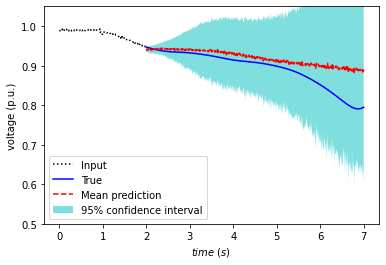}
\end{subfigure}
\caption{\textit{Top:} B-DeepONet mean prediction and confidence interval~(CI) for two post-fault trajectories with inputs contaminated with noise and selected at random from~$\mathcal{D}_\text{test}$. The CI contains $\epsilon_\text{ratio} = 100 \%$ of the true test trajectories. \textit{Bottom:} Prob-DeepONet mean prediction and confidence interval for two post-fault trajectories with inputs contaminated with noise and selected at random from~$\mathcal{D}_\text{test}$. The CI in this noisy scenario contains $\epsilon_\text{ratio} = 100 \%$ of the true test trajectories.}
\label{fig:noisy-confidence-intervals}
\end{figure}

Figure~\ref{fig:noisy-confidence-intervals} depicts the predicted confidence intervals for four test trajectories. The results show that the uncertainty estimate, represented by the $95\%$ confidence interval, does not deteriorate in this noisy scenario.  

We conclude this experiment by showing that B-DeepONet and Prob-DeepONet are robust to noisy inputs. To this end, we report in Table~\ref{table:noisy-ratios} the~$\epsilon_\text{ratio}$ for 100 test trajectories selected from $\mathcal{D}_\text{test}$ whose inputs are contaminated with noise. The results show no deterioration; the $95\%$ confidence intervals still capture about $94 \%$ of the true post-fault trajectories. 

\begin{table}[h!]
\centering
\begin{tabular}{| c | c  c |}
\hline
 & B-DeepONet & Prob-DeepONet \\
\hline
\hline
$\epsilon_\text{ratio}$ & 92.87 \% & 93.27 \% \\
\hline
\end{tabular}
\caption{The average $\epsilon_\text{ratio}$ for 100 post-fault trajectories selected from~$\mathcal{D}_\text{test}$ whose inputs are contaminated with noise.}
\label{table:noisy-ratios}
\end{table}

\subsection{Experiment 4: A power grid application}  \label{sub-sec:alarms}
This final experiment verifies how well the proposed B-DeepONet and Prob-DeepONet predict violations of pre-established under-voltage profiles. Such a violation could trigger, for example, under-voltage load shedding. We refer the interested reader to~\cite{zhang2021adaptive, mozina2007undervoltage} and references therein for more information about these voltage profiles. 

\begin{figure}[t!]
\centering
\begin{subfigure}[b]{0.49\textwidth}
\centering
\includegraphics[width=1.0\textwidth, height=5.5cm]{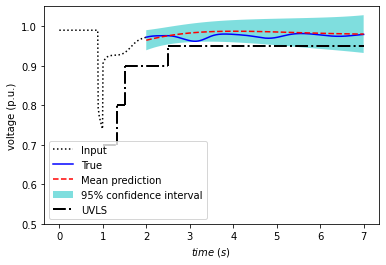}
\end{subfigure}
\begin{subfigure}[b]{0.49\textwidth}
\centering
\includegraphics[width=1.0\textwidth, height=5.5cm]{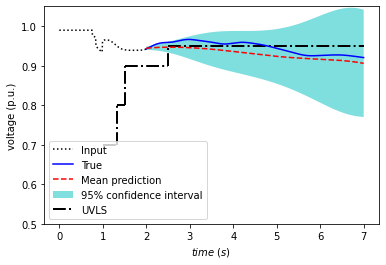}
\end{subfigure}
\caption{The 95\% confidence interval predicted by the B-DeepONet (\textit{left}) and Prob-DeepONet (\textit{right}) compared with the profile for under-voltage load shedding.}
\label{fig:voltage-profile}
\end{figure}

\begin{figure}[t!]
\centering
\begin{subfigure}[b]{0.49\textwidth}
\centering
\includegraphics[width=1.0\textwidth, height=5.5cm]{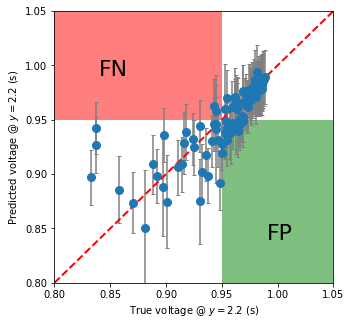}
\end{subfigure}
\begin{subfigure}[b]{0.49\textwidth}
\centering
\includegraphics[width=1.0\textwidth, height=5.5cm]{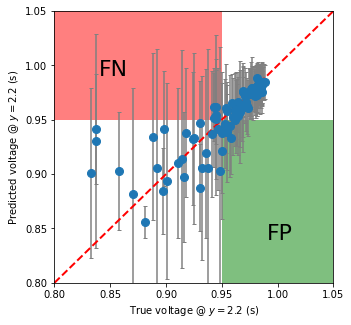}
\end{subfigure}
\caption{A regression scatter plot of the B-DeepONet (\textit{left}) and Prob-DeepONet (\textit{right}) mean predictions (dot) and $95\%$ confidence intervals (error bar) for all test trajectories at time $y = 2.2$ seconds. We verify whether the predicted CIs intersect regions that can trigger under-voltage load shedding. .}
\label{fig:alarms}
\end{figure}

We illustrate in Figure~\ref{fig:voltage-profile} a test trajectory for B-DeepONet and Prob-DeepONet and their respective voltage profiles. One can observe that the $95\%$ confidence interval predicted by B-DeepONet and Prob-DeepOnet may intersect the voltage profiles at some point. Furthermore, we present in Figure~\ref{fig:alarms} a regression scatter plot of the mean predictions and $95\%$ confidence intervals for all test trajectories at time $y = 2.2$ seconds. The figure compares the mean prediction with the true value at time~$y=2.2$ seconds. It also verifies whether the CIs intersect regions that can trigger under-voltage load shedding. In particular, we say a false-negative~(FN) alarm occurs when the whole CI lies within the FN region. For a false-positive~(FP) alarm, we consider two scenarios. 1) In the conservative scenario, we say an FP occurs when the CI intersects the FP region. 2) In a non-conservative scenario, the FP occurs when the whole CI lies within the FP region. 

The results of Figure~\ref{fig:alarms} illustrate that the proposed methods eliminate all false-negative alarms. Moreover, the results show that the methods yield a $32\%$ for B-DeepONet and $29\%$ for Prob-DeepONet (resp. 0\%) false-positive rate in the conservative (resp. non-conservative) scenario. 

\section{Discussion and Future Work} \label{sec:discussion}
We start this section by providing a summary of our results.

\textit{Summary of our results.} Our first main result illustrates (in Figure~\ref{table:predictive-power}) that B-DeepONet and Prob-DeepONet provide comparable mean prediction and generalization than the vanilla DeepONet trained with Adam. Our second main result validates (in Figure~\ref{fig:confidence-intervals} and Table~\ref{table:ratios}) that B-DeepONet and Prob-DeepONet effectively quantify the aleatoric uncertainty, and, hence, they can provide reliable post-fault trajectory predictions. In particular, the predicted confidence intervals capture almost $95\%$ of the true test trajectories. We achieve these results using a small training dataset. Our final result illustrates (in Figure~\ref{fig:voltage-profile} and Figure~\ref{fig:alarms}) that the reliable predictions of the proposed methods enable reducing (to almost zero) the number of false-negative alarms. Such alarms are raised when the methods wrongly predict that a post-fault voltage trajectory violates pre-defined limits. 

Our future work seeks to improve the design of DeepONets for transient stability prediction as follows.

\textit{On relaxing the DeepONet Sensors.} We note that in practice, the set of sensor evaluation points for DeepONets $(x_1,\ldots,x_m)$ (which corresponds to the set of sampling times of a transient response simulator or PMUs) can be affected by delays/noise or not known a-priori. Thus, it is part of our future work to relax the sensor location constraint by designing a resolution-independent DeepONet. 

\textit{On using neighboring information.} In this paper, we did not use trajectory information from neighboring buses to train our DeepONets. In our previous work~\cite{li2020machine}, we observed that using this neighboring information improves training and generalization. Thus, it is part of our future work to design a multi-input multi-output DeepONet.

\textit{On developing a transfer learning strategy.} We plan to design a transfer learning strategy that enables using  a trained B-DeepONet or Prob-DeepONet to train efficiently multiple post-fault trajectory predictors on a large-scale power grid. 

\textit{On power-grid informed DeepONets.} Physics-Informed Neural Networks~(PINN)~\cite{raissi2019physics} are neural networks that we train to satisfy the underlying mathematical physics of the given problem. PINNs have provided us with multiple transformative results by using the physics of the underlying problem as a regularizing agent. This regularization limits the space of possible solutions, which, in turn, enables generalizing well even when the amount of training data is small. In our previous work~\cite{moya2021dae}, we developed a PINN to simulate stiff and nonlinear differential-algebraic equations. In our ongoing work, we are designing DeepONets that use the power grid dynamics as a regularizing agent. We aim to use these power grid-informed DeepONets to enhance generalization capabilities for transient analysis and build a transient analysis tool for large-scale power grids.

\section{Conclusion} \label{sec:conclusion}
This paper proposes using Deep Operator Networks~(DeepONet) to predict post-fault trajectories reliably. In particular, we design a DeepONet to (1) take as inputs the fault-on trajectories collected, for example, via simulation or phasor measurement units, and (2) provide as outputs the predicted post-fault trajectories. Furthermore, we design {\color{black}two} uncertainty quantification methods that provide the much-needed ability to balance efficiency with reliable/trustworthy predictions.  First, we design a \textit{Bayesian DeepONet}~(B-DeepONet) that uses stochastic gradient Hamiltonian Monte-Carlo to sample from the posterior distribution of the DeepONet parameters. Second, we {\color{black}design} a \textit{Probabilistic DeepONet}~(Prob-DeepONet) that uses a probabilistic training strategy to equip DeepONets with a form of automated uncertainty quantification, at virtually no extra computational cost. Finally, we use a series of experiments on the IEEE-{\color{black}68}-bus system that provide the empirical evidence about the effectiveness of the proposed B-DeepONet and Prob-DeepONet.

\section*{Acknowledgment}
The authors gratefully acknowledge the support of the National Science Foundation (DMS-1555072, DMS-1736364, DMS-2053746, and DMS-2134209), and Brookhaven National Laboratory Subcontract 382247, and U.S. Department of Energy (DOE) Office of Science Advanced Scientific Computing Research program DE-SC0021142.

\bibliographystyle{abbrv}
\bibliography{references}

\section*{Appendix 1}
\begin{figure}[t!]
\centering
\begin{subfigure}[b]{0.45\textwidth}
\centering
\includegraphics[width=1.0\textwidth, height=5cm]{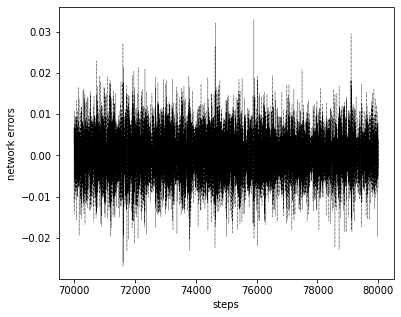}
\end{subfigure}
\begin{subfigure}[b]{0.45\textwidth}
\centering
\includegraphics[width=1.0\textwidth, height=5cm]{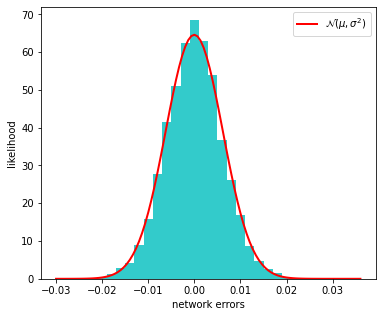}
\end{subfigure}
\\
\begin{subfigure}[b]{0.45\textwidth}
\centering
\includegraphics[width=1.0\textwidth, height=5cm]{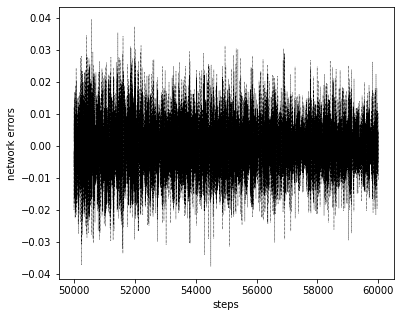}
\end{subfigure}
\begin{subfigure}[b]{0.45\textwidth}
\centering
\includegraphics[width=1.0\textwidth, height=5cm]{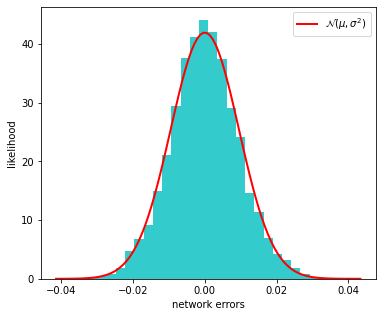}
\end{subfigure}
\\
\begin{subfigure}[b]{0.45\textwidth}
\centering
\includegraphics[width=1.0\textwidth, height=5cm]{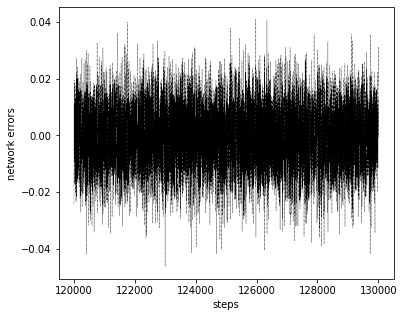}
\end{subfigure}
\begin{subfigure}[b]{0.45\textwidth}
\centering
\includegraphics[width=1.0\textwidth, height=5cm]{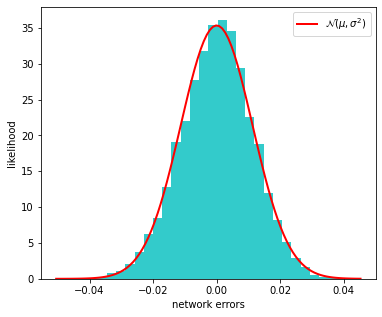}
\end{subfigure}
\caption{\textit{Left:} network errors experienced by the vanilla DeepONet trained with Adam. \textit{Right:} distribution of the network errors. \textit{Top:} network errors and their distribution when $\mathcal{D}_\text{train}$ is simulated using single-line contingencies. \textit{Middle:} network errors and their distribution when $\mathcal{D}_\text{train}$ is simulated using double-line contingencies. \textit{Bottom:} network errors and their distribution when $\mathcal{D}_\text{train}$ is simulated using both single-line and double-line contingencies.}
\label{fig:network-errors}
\end{figure}

In this Appendix, we empirically show that the network errors $(G_\theta(u)(y) - G^\dagger(u)(y))$ the vanilla DeepONet experiences, during gradient descent training (\eg Adam), follows a Gaussian distribution.  To this end, we plot the network errors and their distributions (see Figure~\ref{fig:network-errors}) after some burn-in number (80000) of steps\footnote{One should note that a step is different from an epoch. Here we define \textit{step} as one gradient descent execution.}. In Figure~\ref{fig:network-errors}, we have considered three different cases. In the first case, the training data~$\mathcal{D}_\text{train}$ contains voltage fault-on and post-fault trajectories simulated from single-line contingencies. In the second case, $\mathcal{D}_\text{train}$ contains voltage trajectories simulated from double-line contingencies. In the last case,  $\mathcal{D}_\text{train}$ contains voltage trajectories simulated from both single-line and double-line contingencies. The results show that regardless of the case, the network errors experienced by DeepONets are Gaussian distributed. Such a result provides us with empirical confirmation that the outputs of the Prob-DeepONet can be modeled as the parameters of a Gaussian distribution~$\mathcal{N}(\mu_\theta(u)(y), \sigma_\theta^2(u)(y))$.

\end{document}